\input amstex
 \documentstyle{amsppt}
\loadbold
\topmatter
\pageheight{7.5in}
\title The Fine Structure of the Kasparov Groups II: \\
 topologizing the UCT   \\
      \endtitle
 \rightheadtext{Fine Structure II }
\author    Claude L. Schochet
  \endauthor
 \affil
 Mathematics Department   \\
 Wayne State University    \\
 Detroit, MI 48202      \\ \\
 \endaffil
\address       Wayne State University     \endaddress
\email     {    claude\@math.wayne.edu }     \endemail
\date {May 30, 2001 }\enddate

\endkeywords
  \endsubjclass
\abstract { The Kasparov groups $KK_*(A,B)$ have a natural
structure as pseudo\-polonais groups. In this paper we analyze
how this topology interacts with the terms of the Universal
Coefficient Theorem (UCT) and the splittings of the UCT
constructed by J. Rosenberg and the author, as well as its
canonical  three
term decomposition which exists under bootstrap hypotheses.
 We show that
 the various topologies
on $Ext_{\Bbb Z}^1(K_*(A), K_*(B) )$ and other related groups
mostly coincide.
 Then we
focus attention on the Milnor sequence and the fine structure
subgroup of $KK_*(A,B)$.
An important consequence of our work is that under bootstrap hypotheses
the closure of zero of $KK_*(A,B)$ is  isomorphic to the group
$Pext_{\Bbb Z}^1(K_*(A), K_*(B) )$.
Finally, we introduce new splitting obstructions for the Milnor
and Jensen sequences and prove that these sequences split if
$K_*(A)$ or $K_*(B)$ is torsionfree. }
\endabstract

\endtopmatter

\document
\magnification = 1200
\pageno=1

\def\KK #1.#2 { KK_*( #1,#2) }
 \def\KKgraded #1.#2.#3 { KK_{#1}( #2, #3) }
 \def\bKKgraded #1.#2.#3 {\overline{ KK}_{#1}( #2, #3) }
  \def\dirlim{\underset \longrightarrow\to {lim} \, }
  \def\invlim{\underset \longleftarrow\to {lim} \, }

\def\invlimone{\underset \longleftarrow\to {lim^1}  \, }

\def\ext #1.#2 {Ext _{\Bbb Z}^1( {#1} , {#2} ) }
\def\usualext {Ext _{\Bbb Z}^1 (  K_*(A), K_*(B) ) }
\def\usualpext {Pext _{\Bbb Z}^1 (  K_*(A), K_*(B) ) }
\def\pext #1.#2 {Pext _{\Bbb Z}^1( {#1} , {#2} ) }
\def\hom #1.#2 {Hom_{\Bbb Z} (#1 ,  #2  ) }

 \def\zero {  \overline {  \{   0 \} }  }
 \def\zzero #1{\zero _{#1}}

\newpage\beginsection {Introduction  }

This is the second in a series of papers devoted to the analysis
of the topological structure of the Kasparov groups.
In the first paper in this series \cite {13} we demonstrated
the following facts:

\medbreak
\roster
\item There is a natural structure of a pseudo\-polonais\footnote{
A topological space is {\it {polonais}} if it is separable, complete,
and metric. It is {\it {pseudo\-polonais}} if it is separable and
it has a pseudo\-metric (all axioms for a metric space except
that if $d(x,y) = 0$ then perhaps $x \neq y$) and if its Hausdorff
quotient metric space is polonais.
 If it is a topological group then we insist that the
metric be invariant.}
group
on $\KKgraded *.A.B  $.  \cite {13, 2.3}
\medbreak
\item The Kasparov pairing is jointly continuous with respect
to this topology. \cite {13, 3.8}
\medbreak
\item The index map
$$
\gamma : \KKgraded *.A.B  \to \hom {K_*(A)}.{K_*(B)}
$$
is continuous.
If $Im(\gamma ) $ is closed (e.g., if $\gamma $ is onto), then
$\gamma $ is an open map.
 If $\gamma $ is an algebraic isomorphism then it is an
isomorphism of topological groups. \cite {13, 4.6}

\endroster

\medbreak

In this paper we will study the various purely algebraic and
also analytic topologies which naturally occur on the components
of $\KKgraded *.A.B $ within the context of the UCT. Special
attention will be paid to $Z_*(A,B)$, the closure of zero
in $\KKgraded *.A.B $, which we call the fine structure subgroup.

\medbreak

Section 1 briefly summarizes various results on the structure
 of $\KKgraded *.A.B $ provided that $A \in \Cal N$, the bootstrap
 category. The canonical KK-filtration diagram (1.4) is introduced.
 This diagram is used heavily in the remainder of the paper.

Section 2 deals with several \lq\lq algebraic\rq\rq\,   topologies
on $\ext G.H $ where $G$ and $H$ are (usually countable) abelian
 groups. We show that the Jensen isomorphism
$$
\pext G.H \,\cong\, \invlimone \hom {G_i}.H
$$
(which holds for $G$ written as an increasing union of finitely
generated subgroups $G_i$) is an isomorphism of topological
groups.

 Section 3 deals with topologizing most of the terms in the
KK-filtration diagram. The key result is
the natural isomorphism
$$
Z_*(A,B)  \,\cong\,   \pext {K_*(A)}.{K_*(B)}
 $$
where  $Z_*(A,B) $ denotes
the closure of zero in $\KKgraded *.A.B $.

In Section 4 we complete our analysis of the topological
structure of the KK-filtration diagram. We show that each of
the algebraic splittings of the UCT sequence is continuous.

In Section 5 we introduce a new type of invariant  for
pairs of $C^*$-algebras $A$ and $B$, namely, splitting
invariants
$$
m(A,B) \in Ext_{\Bbb Z}^1(   \,
\invlim \KKgraded *.{A_i}.B   \, ,  \,  \invlimone
\KKgraded *.{A_i}.B\, )
$$
and
 $$
j(A,B) \in
Ext_{\Bbb Z}^1 (\, \invlim\ext{K_*(A_i)}.{K_*(B)} \,  ,  \,
\pext K_*(A).K_*(B) \,  ).
$$
The obstruction $m(A,B)$ vanishes iff the Milnor sequence
splits, and the obstruction $j(A,B)$ vanishes iff the Jensen sequence
 splits. We prove that
$$
j(A,B) =(\invlim \delta _i)^* m(A,B)
$$
up to isomorphism
so that the splitting of the Milnor sequence implies the splitting
of the Jensen sequence.  If $K_*(A)$ or $K_*(B)$ is torsionfree then
both invariants are zero and hence both sequences split. To find
non-splitting examples one must assume that both  $K_*(A)$ and
 $K_*(B)$ have $p$-torsion for some prime $p$. We briefly examine
this situation.

\medbreak

In this paper all $C^*$-algebras are separable.
 All $C^*$-algebras appearing in the first
variable of $KK$ are assumed to be nuclear.  Whenever
$\ext {K_*(A)}.{K_*(B)} $ is taken to be a subgroup of $\KKgraded *.A.B $
it is understood that $A \in \Cal N$ and that the inclusion is
via the UCT. All pseudo\-polonais and polonais groups
are understood to be abelian.

\medbreak

Note: Marius Dadarlat  recently has discovered \cite{4, Cor. 4.6}
  a different
and very interesting proof of
Theorem 3.3 as well as some of its consequences for
quasi\-diagonality which we first proved in this paper
and in \cite{14}.  It is a pleasure to thank him as well
as his former student Nathaniel Brown for helpful advice and
assistance in this series of papers.


\newpage\beginsection{ 1. Algebraic Facts}

In this section we briefly summarize various results concerning
 $\KKgraded *.A.B $.
We assume that $A$ is in the bootstrap category $\Cal N$ so
that the UCT \cite {9} holds. We introduce the
canonical KK-filtration diagram (1.4) which is central
to our analysis.

The first result is purely algebraic.
We recall that for abelian groups $G$ and $H$, $\pext G.H $
denotes the subgroup of $\ext G.H $ consisting of pure
extensions.\footnote{An extension of abelian groups
$$
0 \to H \to E  \to G \to 0
$$
is {\it{pure}} if $H \cap nE = nH $ for each positive integer $n$.
Equivalently, the extension
is pure if it
 splits when restricted to every
finitely generated subgroup of $G$. For further information on $Pext$
please see \lq\lq A Pext Primer" \cite {16}.}

\proclaim{Theorem 1.1, C. U. Jensen \cite{6}}
Let  $G$ be a countable group written as a
union of
a sequence of
finitely generated subgroups $G_i$ and let $H$ be any group.
Then there is a natural isomorphism
$$
\pext G.H   \cong  \invlimone \hom {G_i}.H   .
$$
\endproclaim\qed

We  recall the definition of a
$KK$-filtration from \cite {11}.
\medbreak

\proclaim {Definition 1.2}  A $KK$-filtration
 of a separable $C^*$-algebra
$A$ is an increasing sequence of commutative $C^*$-algebras
$$
A_0 \hookrightarrow A_1 \hookrightarrow A_2 \hookrightarrow \dots
$$
which satisfies the following conditions:
\roster
\item  $A_i^+ \cong C(X_i)$ for some finite $CW$-complex $X_i$.
\medbreak
\item Each map $K_*(A_i) \to K_*(A_{i+1})$ is an inclusion.
\medbreak
\item $\dirlim A_i$ is $KK$-equivalent to $A$.
\medbreak
\endroster
\endproclaim
\medbreak

It follows that each $K_*(A_i)$ is finitely generated and that
$$
\dirlim K_*(A_i) \,\cong\, K_*(A)
$$
so that the sequence $\{ K_*(A_i)  \} $ is an increasing sequence of
finitely generated subgroups with limit $K_*(A)$.
Since the UCT is preserved under $KK$-equivalence, it follows that
any $KK$-filtered $C^*$-algebra $A$ satisfies the UCT for all $B$.
We show \cite {11,  1.5} that each $A \in \Cal N$ has a
$KK$-filtration and that the filtration is unique in the sense that
groups such as
$$
\invlimone \hom {K_*(A_i)}.{K_*(B)}
$$
which {\it {a priori }} depend upon  a choice of $KK$-filtration,
in fact depend only upon $K_*(A)$ and are independent of choice
of $KK$-filtration.
\medbreak

\medbreak

 \proclaim { Theorem  1.3  \cite {11, Theorem 1.6} } Suppose
that $A$ has $KK$-filtration $\{ A_i \}$.
   Then the following
  diagram
\medbreak
$$
\minCDarrowwidth {.155in}
\CD
@. @.    @.      0   @.          \\
@.    @.       @.      @VVV            @.    \\
@.       0          @.    @.   \invlim\ext {K_*(A_i)}.{K_*(B)}     @.  \\
@.       @VVV                 @.        @VV{\invlim\delta _i}V
@.         \\
 0  @>>>  \invlimone\KKgraded *.{A_i}.B  @>{\sigma}>>
 \KKgraded *.A.B  @>{\rho}>>   \invlim\KKgraded *.{A_i}.B  @>>>  0  \\
@.   @VV{\psi}V         @VV{id}V     @VV{\tilde\gamma}V   @.   \\
 0  @>>>    \ext {K_*(A )}.{K_*(B)}
 @>{\delta}>> \KKgraded *.A.B
 @>{\gamma }>>  Hom_{\Bbb Z}(K_*(A), K_*(B))  @>>>   0   \\
@.  @VV{\varphi}V   @.    @VVV   @.   \\
@.     \invlim\ext {K_*(A_i)}.{K_*(B)}    @.  @.  0    @.    \\
@. @VVV    @.  @.  @.   \\
@.  0   @.  @.  @.
\endCD
\tag 1.4
$$
\medbreak
 \flushpar is commutative, is natural with respect to $A$ and $B$
and  has exact rows and columns. Each of the groups is independent
of choice of $KK$-filtration and depends only upon $K_*(A)$ and $K_*(B)$.
 Further, the UCT map $\gamma $
splits  unnaturally .
\endproclaim
\medbreak

Note that the lower row is simply the UCT for $\KKgraded *.A.B $
and the row above it is the Milnor sequence associated to the
given $KK$-filtration.
The right column arises as the inverse limit of UCT sequences
for $(A_i,B)$. The map $\tilde\gamma $ is onto since
by a result of Roos,
(cf. \cite{8})
$$
\invlimone \ext{G_i}.H = 0.
$$
The exactness of the left column is also a consequence of a theorem of
Roos \cite{8}.

An immediate consequence of the UCT is the following corollary.

\medbreak
\proclaim{Proposition 1.5} If $K_*(A)$ is finitely generated
then the group $\KKgraded *.A.B $ is countable.
\endproclaim

\medbreak
\demo{Proof}
In light of the UCT,
it suffices to demonstrate that
if $G$ is finitely generated and $H$ is countable then
 both $\hom G.H $
and $\ext G.H $ are countable, and these are elementary.
\enddemo\qed

\newpage


\beginsection{ 2. Algebraic topologies on $Ext$.  }

In this section we introduce several topologies on $\ext G.H $
and determine how they are related. We introduce
a topology on $\invlimone $ and show that the Jensen isomorphism
is an isomorphism of topological groups.

\medbreak

\proclaim {Definition 2.1}
 Suppose that $G$ and $H$ are countable abelian groups.
 We consider three natural topologies on the group $\ext G.H $.
They are:
\medbreak
\roster
\item The $\Bbb Z$-adic topology, where the subgroups $n\ext G.H $
are taken as a system of neighborhoods of the identity.
This will be denoted $\ext G.H  _{\Bbb Z} $.
\medbreak

\item The\footnote{Actually there is one topology for every choice of
injective resolution, but it will be clear from the proof that
these all yield homeomorphic topologies. }
quotient topology obtained from regarding $Ext $ as
a quotient group of a   $Hom $ group resulting from an injective resolution of
$H$.
This will be denoted $\ext G.H _I  $.
\medbreak

\item The Jensen topology, taking the subgroups
$$
Ker {\biggl[} \ext G.H     \longrightarrow   \ext {G_i}.H  {\biggl ]}
$$
as neighborhoods of the identity. This will be denoted $\ext G.H _J $.
\endroster
\endproclaim
\medbreak
 The second topology will be defined precisely as it arises.

We begin by some elementary observations.
If $G$ is a pseudo\-polonais group,
let $G_o$ denote the closure of zero in $G$
and $\overline G = G/G_o $ denote the quotient group.
Note that any algebraic isomorphism $G_o \cong H_o$ is an
isomorphism of topological groups, since these groups are
indiscrete.

We topologize $\hom G.H $ by the topology of pointwise convergence.

\medbreak

\proclaim{Proposition 2.2}
\roster
\item $\hom {-}.{-} $ is a bifunctor from discrete
abelian
 groups to polonais
groups and continuous  homomorphisms in each variable.
\medbreak
\item $\ext {-}.{-} _I $ is a bifunctor from dis\-crete
abelian  groups to
pseudo\-polonais
groups and continuous  homomorphisms.
\medbreak
\item The boundary homomorphisms in each variable
in the respective
 $Hom$-$Ext _I$ sequences are
continuous.
\medbreak
\item If $G$ is finitely generated then $\ext G.H  $ is discrete
in   the $I$ and $J$ topologies.

\endroster
\endproclaim
\medbreak

\demo{Proof}  Part 1) is immediate.

For Part 2), let
$$
0 \to H \to I \to I' \to 0
$$
be an injective resolution of $H$. Then $\ext G.H  _I$ is
the quotient of the polonais group $\hom G.{I'} $ by the
image of the polonais group $\hom G.I $. The image may
not be a closed subgroup, and hence the quotient $\ext G.H $
is pseudo\-polonais but not necessarily polonais.

 Suppose that $r: H_1 \to  H_2 $. Choose  an injective resolution
for $H_1$ and another for $ H_2$ so there is a commuting diagram
 $$
\CD
0 @>>>   H_1    @>>>   I_1   @>>>   I_1'   @>>>   0 \\
@.  @VV{r}V   @VV{s}V     @VV{t}V    \\
0 @>>>  H_2  @>>>  I_2   @>>>    I_2'   @>>>   0
\endCD
$$
(The maps $s$ and $t$ exist since $I_2$ is injective). This induces
a commuting diagram
$$
\CD
\hom G.{I_1'}   @>{t_*}>>     \hom G.{ I_2'}   \\
@V{p_1}VV    @V{p_2}VV    \\
\ext G.{H_1}  _I  @>{r_*}>>   \ext G.{H_2}   _I
\endCD
$$
The maps $t_*$, $p_1$ and $p_2$ are continuous,  and
 $p_1$ and $p_2$ are quotient maps. This implies that $r_*$ is
continuous.

 For Part 3)  we argue as
follows. Suppose that $G$ is fixed and
$$
0 \to H'  \to H \to H'' \to 0
$$
is a short exact sequence. We are asked to prove that the boundary
homomorphism
$$
\delta ^H :  \hom G.{H''} \to \ext G.{H'}
$$
is continuous. Let
$$
0 \to H' \to I \to I'  \to 0
$$
be an injective resolution of $H'$. Then there is a commuting
diagram
$$
\CD
0 @>>>   H'   @>>>   H   @>>>   H''   @>>>   0 \\
@.  @VV{id}V   @VV{f}V     @VV{g}V    \\
0 @>>>   H'   @>>>   I   @>>>   I'   @>>>   0
\endCD
$$
(this exists since $I$ is injective) and by the naturality of the
long exact sequence there is a commutative square
$$
\CD
\hom G.{H''}    @>{\delta ^H}>>    \ext G.{H'}  _I \\
@VV{g_*}V      @VV{id}V   \\
\hom G.{I'}    @>{\delta ^I}>>   \ext G.{H'}  _I
\endCD
$$
so that $\delta ^H = \delta ^Ig_*$. The map $g_*$  is continuous
 by part 1) and the map $\delta ^I $ is continuous since
it is the map which
is used to define
 the quotient topology to $\ext G.H   _I$.
So $\delta ^H $ is continuous  as required.
The continuity of the functor $\ext {-}.H   _I $ is immediate
from definitions.

Part 4) holds since the group $\ext G.H $ is countable
and complete in each of the topologies, hence discrete.
\qed\enddemo

\medbreak

\proclaim {Theorem 2.3}
 For any abelian groups $G$ and $H$,
$$
\overline {  \{ 0  \} }_I \,=\,
\overline {  \{ 0  \} }_J \,=\,
\overline {  \{ 0  \} }_{\Bbb Z}  \,
\,=\, \pext G.H \,
\subseteq
    \, \ext G.H   .
$$
\medbreak
\medbreak

\endproclaim
\demo {Proof} The identity
$$
\zzero {\Bbb Z}  \cong \pext G.H
$$
is more or less by definition. Roos's Theorem \cite{8} gives
us the isomorphism
$$
\zzero J   \cong \invlimone \hom {G_i}.H
$$
and then applying Jensen's Theorem (1.1) yields the identity
$$
\pext G.H \cong \zzero J  .
$$

It remains to identify $\zzero I $.
Suppose that $a \in \zzero I. $
Then
for each $i$,
$$
\varphi _i(a) = 0 \in \ext{G_i}.H _I    \, \cong\, \ext{G_i}.H  _J
$$
since this group is Hausdorff, and hence
$$
a \in Ker \bigg[ \varphi :  \ext G.H
 \longrightarrow \invlim \ext{G_i}.H  \bigg]  \cong \zzero J .
$$
Thus $\zzero I \subseteq \zzero J $.

In the other direction, suppose that $a \in \zzero J $.
Let
$$
0 \to G \to I \overset\zeta\to\longrightarrow I'  \to 0
$$
be an injective resolution of $G$.
Then $\varphi (a) = 0$, and hence $a|_{G_i} = $ for each $i$.
Represent $a = [f]$ for some $f: G \to I'$. Then
$$
[f|_{G_i} ] = 0 \in \ext{G_i}.H
$$
for each $i$. Thus there exists functions $h_i : G_i \to I$
such that the diagram
$$
\CD
G_i @>>>   G  \\
@VV{h_i}V  @VV{f}V  \\
I @>{\zeta}>>   I'
\endCD
$$
commutes. Since $I$ is injective we may extend the map $h_i$
to a map $\hat h_i : G \to I$. Then $\{ \hat h_i  \}$ is a
sequence in $\hom G.I $.

We claim that $\hat h_i \to f$ in the topology of pointwise
convergence. This is easy. We must show that for each
$x \in G$ that $\hat h_i(x)  \to f(x)$. This is true since
as soon as $i$ is large enough so that $x \in G_i$ we have
$\hat h_i(x)  = f(x)$.

This proves that $\zzero J \subseteq \zzero I $. Combining with
the first part of the proof yields $\zzero J = \zzero I $
which completes the proof of the theorem.
\qed\enddemo

\medbreak

{\bf {Remark 2.4}.} In general the topological groups  $\ext G.H _{\Bbb Z} $
and $\ext G.H _{S} $ are not homeomorphic. Here is an example, courtesy of
C. U. Jensen. Let $G_i $ be the direct sum of $i$ copies of the group
$\Bbb Z/2 $  and $G_i \to G_{i + 1}$ be  the map $x \to (x,0)$. Then
$G = \dirlim G_i $ is the direct sum of countably many copies of
$\Bbb Z/2$, and so
$$
\ext  G.{\Bbb Z}
 \,\cong\,  \prod _{i = 1}^{\infty }
\ext  {\Bbb Z/2}.{\Bbb Z} \,\cong\,
 \prod _{i = 1}^{\infty } {\Bbb Z/2}    .
$$
The $\Bbb Z$-adic topology is discrete in this case, since
$2\ext G.H = 0$. However,
$$
Ker \, \bigg[   \ext G.{\Bbb Z}  \to \ext {G_n}.{\Bbb Z} \bigg]
\,\cong\,   \prod _{i = n+ 1}^{\infty } \Bbb Z/2
$$
and so $\ext G.H _{S} $ is not discrete.

\medbreak

To complete our algebraic discussion, we consider the groups
that arise from writing $G$ as a union of an increasing family
of finitely generated subgroups $G_i$.
 The induced inverse sequence $\{\,\hom {G_i}.H \,\}$
has associated to it the canonical Eilenberg exact sequence
 which we may use to define $\invlimone $.
\medbreak
$$
\minCDarrowwidth {.16in}
\CD
0 @>>>   \hom {G}.H     @>>>  \prod _i \hom {G_i}.H
@>\Psi >>   \\ @. \\
@.
@>\Psi >>  \prod _i \hom {G_i}.H
@>>>
\invlimone \hom {G_i}.H  @>>> 0 .
\endCD
\tag 2.5
$$
This definition is independent of choice of subgroups.
We give $\prod _i \hom {G_i}.H $ the product topology and
give
$\invlimone \hom {G_i}.H $ the  quotient topology.

\medbreak

\proclaim{Proposition 2.6} The Jensen isomorphism (1.1) is an
isomorphism of topological groups
$$
\pext G.H \overset\cong\to\longrightarrow     \invlimone \hom {G_i}.H   .
$$
\endproclaim
\medbreak

\demo{Proof}
We are given a direct sequence of abelian groups
$$
G_0 \to G_1 \to G_2 \to \dots
$$
and this yields the canonical pure short exact sequence
$$
0 \to \oplus G_i \overset\psi\to\longrightarrow
\oplus G_i  \to G \to 0.
$$
  Apply the functor $\hom {-}.H $ and one obtains the following
  commutative diagram with exact columns:

$$
\CD
0  @.   0  \\
@VVV   @VVV   \\
\hom G.H  @>{\cong}>> \invlim \hom{ G_i}.H    \\
@VVV    @VVV    \\
\hom{\oplus G_i}.H   @>{H}>>  \prod \hom{ G_i}.H  \\
@VV{\psi ^*}V    @VV{\Psi }V    \\
\hom{\oplus G_i}.H   @>{H}>>  \prod \hom{ G_i}.H  \\
@VVV    @VVV   \\
\pext G.H  @>{h}>>  \invlimone \hom{ G_i}.H   \\
@.   @VVV   \\
@.    0
\endCD
$$
where $\Psi $ is the canonical Eilenberg map
and $H$ is the evident isomorphism of topological groups. Then
the maps $h$  and $h^{-1}$ are  continuous bijections, hence
isomorphisms of
topological groups.
 \enddemo\qed


\newpage
\beginsection{3. Topologizing the $KK$-filtration diagram}

In this section we topologize most of the  terms in the $KK$-filtration
diagram (1.4), show that each of the natural maps is continuous,
that some are open, and as a consequence obtain a deeper
understanding of the fine structure subgroup. All maps in
this section not otherwise identified are identified
in the KK-filtration diagram (1.4). The most important result
for applications is Theorem 3.3, where we identify $Z_*(A,B)$,
the closure of zero in the Kasparov groups, as
$\pext {K_*(A)}.{K_*(B)} $.

We begin by topologizing $\KKgraded *.A.B $ with its natural analytic
topology, as described in detail in \cite{13}. (This arises
from giving extensions the topology of pointwise convergence.)
With respect to this topology $\KKgraded *.A.B $ is a
pseudo\-polonais topological group.

Let $Z_*(A,B)$ denote the closure of zero in $\KKgraded *.A.B $.
Any  KK-equivalence of $A$ to $A'$ induces an algebraic isomorphism
$$
\KKgraded *.{A'}.B  \,\cong\,  \KKgraded *.A.B   .
$$
Since KK-pairings are continuous (this is the principal result
of \cite {13}),  it follows \cite{13, Theorem 5.1} that a KK-equivalence
induces an isomorphism of topological groups. The closure of
zero is of course preserved under such an isomorphism, and
hence
$$
Z_*(A,B) \,\cong\, Z_*(A', B)  .
$$
The analogous result holds in the second
variable.
\medbreak

Every $C^*$-algebra in the bootstrap
category $\Cal N$ has a $KK$-filtration which is unique up to
$KK$-equivalence. A $KK$-equivalence induces an isomorphism
of topological groups, by
\cite{13, 5.1} and hence up to isomorphism of topological
groups we may assume
without loss of generality that each $A$ has a $KK$-filtration,
and the entire KK-filtration diagram depends only upon
$A$ and $B$ and is natural in each variable.

The group $\hom {K_*(A)}.{K_*(B)} $ is topologized
as before
by first
declaring $K_*(A)$ and $K_*(B)$ to be discrete and then by using the
topology of pointwise convergence on $Hom$. We may regard it as a
subgroup of the countable (one for each element of $K_*(A)$)
product of copies of the (discrete) group $K_*(B)$.  Thus $\hom
{K_*(A)}.{K_*(B)} $ has the structure of a polonais group.

We recall from \cite{13, 7.4} that the natural index map
$$
\gamma :  \KKgraded *.A.B \longrightarrow   \hom {K_*(A)}.{K_*(B)}
$$
is continuous. Thus
$$
ker(\gamma ) = \ext {K_*(A)}.{K_*(B)} _{rel}
$$
is a closed subgroup of $\KKgraded *.A.B $ with the relative topology.

\medbreak

Next we topologize the group $\invlim \KKgraded *.{A_i}.B $ by
giving it the relative topology with respect to the natural
inclusion
$$
\invlim \KKgraded *.{A_i}.B    \longrightarrow \prod \KKgraded
*.{A_i}.B   .
$$
 Note that each $\KKgraded *.{A_i}.{B} $ is countable,
complete, and hence discrete. Then it is easy to show that the
maps $\rho $ and $\bar\gamma $ are continuous. Further, the map
$\bar\gamma $ is open since $\gamma $ is open.

There are two reasonable topologies on
 the group $\invlim \ext
{K_*(A_i)}.{K_*(B)} $.
One possibility is to give it
 the relative topology with
respect to the inclusion $\invlim \delta _i$, or, equivalently, by
giving it the relative topology in the group
$$
\prod _i   \ext {K_*(A_i)}.{K_*(B)}   .
$$
This topology we denote by  $\invlim \ext {K_*(A_i)}.{K_*(B)}
_{rel}$.

 Alternately, we may topologize  $\invlim\ext {K_*(A_i)}.{K_*(B)}
$
as in Section 2
by giving it the quotient topology as a quotient of the group $
\ext {K_*(A )}.{K_*(B)} _I $ via the map $\varphi $.
 Denote this option by
 $$
\invlim\ext {K_*(A_i)}.{K_*(B)}  _{I}   .
$$
 Under this option
the map $\varphi $ is continuous and open.

\medbreak

\proclaim {Proposition 3.1} The natural map
$$
\invlim\ext{K_*(A_i)}.{K_*(B)}  _I  \longrightarrow
\invlim\ext{K_*(A_i)}.{K_*(B)} _{rel}
$$
is a homeomorphism.
\endproclaim
\medbreak

\demo {Proof} We refer to diagram (1.4)
  for notation. The map $\rho\delta $
is continuous (essentially by definition, in both cases) and the
composition $\tilde\gamma\rho\delta = 0$, so that the image of
$\rho\delta $ lies in the image of the map $\invlim \delta _i $.
This implies that the natural map
$$
\rho\delta : \ext {K_*(A )}.{K_*(B)} _I \longrightarrow
\invlim\ext {K_*(A_i)}.{K_*(B)}  _{rel}
$$
is continuous. It vanishes, of course, on the image of $\psi $ and
hence produces a continuous map $I$ (the identity!)
$$
I :  \invlim\ext {K_*(A_i)}.{K_*(B)}  _{I} \longrightarrow
\invlim\ext {K_*(A_i)}.{K_*(B)}  _{rel}
$$
which is obviously an algebraic isomorphism.

We may explicitly construct the inverse to $I$ as follows.
 Given any
commuting diagram of abelian groups with exact rows
$$
\CD
0 @>>> G'  @>>>  G  @>>> G''  @>>>  0  \\
@.       @VV\gamma 'V     @VV\gamma V    @VV\gamma V''    @.  \\
0 @>>> H' @>>>  H @>>>   H'' @>>>  0
\endCD
$$
the Snake Lemma asserts that there is a long exact sequence
$$
0 \to  Ker(\gamma ') \to Ker(\gamma ) \to Ker(\gamma '')
\overset\delta\to\rightarrow Cok(\gamma ') \to Cok(\gamma ) \to
Cok(\gamma '') \to 0
$$
We show \cite {15} that in the context of topological groups the
map $\delta $ is a continuous algebraic isomorphism
$$
 \delta :  \invlim\ext {K_*(A_i)}.{K_*(B)}  _{rel}
 \longrightarrow
  \invlim\ext {K_*(A_i)}.{K_*(B)}  _{I}   .
$$
Since $\delta ^{-1} = I$
 our results imply that $I$ is an open map and thus
 is an isomorphism of topological groups.
\qed\enddemo

Henceforth we shall use the notation $\invlim\ext
{K_*(A_i)}.{K_*(B)} $ without subscript to indicate the
topological group.

\medbreak
Thus we conclude:

\proclaim{Proposition 3.2} The right column of the $KK$-diagram is
a short exact sequence of topological groups. The map $\bar\gamma
$ is continuous and open, hence a quotient map, and
$
\invlim \ext {K_*(A_i)}.{K_*(B)}  $ is a closed subgroup of
the group $\invlim \KKgraded *.{ A_i }.B $.
\endproclaim \qed

\medbreak

\proclaim{Theorem 3.3} Suppose that $A \in \Cal N$. Then
the closure of zero in $\KKgraded *.A.B $ is the
group $\usualpext $. That is,
$$
Z_*(A,B)  \,=\, \usualpext    .
$$
\endproclaim

\medbreak
\demo{Proof} In light of Theorem 2.3, it suffices to
show that $Z_*(A,B)  = \zzero I$.

First we show that the identity map
$$
\usualext _I \longrightarrow \usualext _{rel}
$$
is continuous. For this we must recall parts of the
proof of the UCT.

Recall that a geometric injective resolution of $K_*(B)$ is a
sequence
$$
0 \to I_0 \to I_1 \to SB  \to 0
$$
where
   $ I_0$ is an ideal in $I_1$, $K_*(I_j)$ is
injective (that is, divisible) for each $j$,
   and the resulting $K$-theory
long exact sequence degenerates into  an injective resolution
of $K_*(B)$ of the form
$$
0 \to  K_*(B) \to K_*(I_0) \to K_*(I_1) \to 0    .
$$
These exist and are key to our proof of the UCT \cite {9}.

 Consider the natural commuting diagram
$$
\CD
@.    0   \\
@VVV      @VVV   \\
\KKgraded *.A.B    @>\gamma >>    \hom {K_*(A)}.{K_*(B)}   \\
@VV{\hat f_*}V       @VV{f_*}V   \\
 \KKgraded *.A.{I_0}    @>{\gamma _0 }>>    \hom {K_*(A)}.{K_*(I_0)}   \\
@VV{\hat g_*}V       @VV{g_*}V   \\
 \KKgraded *.A.{I_1}    @>{\gamma _1 }>>    \hom {K_*(A)}.{K_*(I_1)}   \\
@VV{\hat\delta _* }V       @VV{\delta _*}V   \\
 \KKgraded *.A.B    @.  \ext {K_*(A)}.{K_*(B)}   \\
@VV{\hat f_*}V      @VVV  \\
@.       0
\endCD
\tag 3.4
$$
The maps $\gamma _1 $ and $\gamma _0$ are isomorphisms  by the UCT,
since the $K_*(I_j)$ are injective,
and thus isomorphisms of topological groups, by \cite {13, 4.6}.
The left column is exact by exactness properties of the Kasparov
groups, and the right column is exact by the usual $Hom$-$Ext$
exact sequence.

Both of the topologies of interest to us arise in this diagram.
Specifically,
$$
 \ext {K_*(A)}.{K_*(B)}  _{rel}   \,\cong\, Im(\hat\delta _*)_{rel}
$$
by the proof of the UCT \cite {9}, and
$$
 \ext {K_*(A)}.{K_*(B)}   _{I}   \,\cong\, Im(\delta _*)_{quot}
$$
by definition.

 It is clear from the diagram that the identity
map
$$
\usualext _I  \longrightarrow  \usualext _{rel}
$$
is continuous. This immediately implies that $\zzero I \subseteq
\zzero {rel} $.

Now consider the canonical map
$$
\varphi : \usualext \longrightarrow
\invlim \ext{K_*(A_i}.{K_*(B)}  .
$$
This map is continuous in both topologies, and the two
induced topologies coincide on  $\ext{K_*(A_i}.{K_*(B)} $
by Proposition 3.1. Thus the map $\varphi _{rel}$
is a continuous homomorphism
from $\usualext _{rel} $ to a Hausdorff topological group. So
$Z_*A,B) $ must be contained in $Ker(\varphi ) $.
We know that
$$
Ker(\varphi ) \cong \usualpext \cong \zzero I
$$
and hence $Z_*(A,B)  \subseteq \zzero I $. Putting
together the information from the previous
paragraph yields $Z_*(A,B)  = \zzero {I} $ as claimed.
\enddemo\qed

Next we  look carefully at the first row of the $KK$-diagram,  the
Milnor $\invlimone $ sequence. As usual we assume that $A \in \Cal N$.

\medbreak

\medbreak \proclaim{Proposition 3.5} The natural map
$$
\rho : \KKgraded *.A.B  \to  \invlim \KKgraded *.{A_i}.B
$$
is continuous and open.
\endproclaim
\medbreak

\demo{Proof}
Continuity is obvious. We must show that $\rho $ is open.
Note that
$$
Ker(\rho ) \cong Z_*(A,B)
$$
by Theorem 3.3.
The proof then comes down to an easy lemma.
\enddemo

\proclaim {Lemma 3.6}
Suppose that $G$ is a pseudopolonais group
with Hausdorff quotient group $\overline G$.  Then the natural map
$\rho : G \to \overline G $
is open.
\endproclaim
\medbreak
\demo{Proof} Suppose that $U$ is an open neighborhood of $0 \in G$. It suffices
to show that $\rho ^{-1}\rho (U) $ is open in $G$. But
$\rho ^{-1}\rho (U) = U$.
\enddemo\qed
\medbreak

\proclaim{Proposition 3.7} The right column of diagram 1.4 is
a short exact sequence of polonais groups. The group
$Im(\invlim \delta _i )$ is a closed subgroup. The function
$\invlim \delta _i $ is an open map onto its image. If $\invlim \delta _i $
is a bijection then it is an isomorphism of topological groups.
\endproclaim
\medbreak

\demo{Proof} Each of the groups in the sequence is polonais
by earlier results. Then
$$
Im(\invlim \delta _i ) = Ker(\bar\gamma )
$$
 is a closed subgroup.
The map $\invlim \delta _i $ is an open map onto its image
by  \cite {13, 6.4(1)}, and if it is a bijection then it
is continuous, open, and hence a homeomorphism.  \qed\enddemo

\newpage
\beginsection{4. Topologizing the UCT}

In this section we complete our topological results and
show how these results fit together in the UCT. We show
that each of the algebraic splittings of the UCT
produced by \cite{9} is continuous. Finally, we demonstrate
 how to realize Jensen's counterexample (2.4) geometrically.

 Our first
task is to study the various $\invlimone $ terms that arise.

The first topology that we consider on the group
$\invlimone \KKgraded *.{ A_i}.B  $ is the relative
topology as a subspace of $\KKgraded *.A.B $ via the map
$\sigma $. We denote this by
$\invlimone \KKgraded *.{ A_i}.B _{rel} $.
It is clear that the Jensen isomorphism
induces a natural isomorphism of topological groups
$$
\invlimone \KKgraded *.{ A_i}.B _{rel}
\,\cong\,
\pext{K_*(A)}.{K_*(B)} _{rel}
$$
and we know by (3.3) that
$$
\pext{K_*(A)}.{K_*(B)} _{rel} \,\cong\, Z_*(A,B)  .
$$

\medbreak
We summarize:
\medbreak

\proclaim {Proposition 4.1} The Milnor sequence
$$
0 \to \invlimone \KKgraded *.{ A_i}.B  _{rel}
\overset\sigma\to\longrightarrow  \KKgraded *.A.B
\overset\rho\to\longrightarrow
 \invlim \KKgraded *.{A_i}.B  \to 0
$$
is a short exact sequence of topological groups, the map
$\rho $ is open and hence a quotient map, and the group
$
\invlimone \KKgraded *.{ A_i}.B \cong Z_*(A,B)
$
is the closure of zero in the analytic topology on $KK$.
\endproclaim

\medbreak

We note that the group    $\invlim \KKgraded *.{A_i}.B  $
appears in the work of R{\o}rdam,  Dadarlat and Loring and is there
denoted $KL_*(A,B)$.   Proposition 4.1  implies that this
group is the Hausdorff quotient of $\KKgraded *.A.B $:
$$
KL_*(A,B) \cong \frac{\KKgraded *.A.B  }{Z_*(A,B) } .
$$
\medbreak

The group $\invlimone \KKgraded *.{ A_i}.B  $ may
also have the quotient topology obtained as a quotient
of the group
$$
\prod \KKgraded *.{ A_i}.B
$$
as in diagram (2.5), but this topology coincides with relative
topology by Proposition 2.6. Similarly there is a natural
isomorphism of topological groups
$$
\invlimone \hom {K_*(A_i)}.{K_*(B)}        _{quot}
\,\cong\,
 \invlimone \hom {K_*(A_i)}.{K_*(B)}      _{rel}
$$
and so we write this group henceforth without subscript.

\medbreak

\proclaim {Theorem  4.2}  Suppose that $\{A_i \}$ is a $KK$-filtration
of $A \in \Cal N$. Then the natural isomorphism
\medbreak
$$
\invlimone\gamma _i : \invlimone \,\KKgraded *.{ A_i }.{ B  }
\longrightarrow
  \invlimone \,\hom {K_*(A_i)}.{K_*(B)}
$$
is an isomorphism of topological groups.

\endproclaim
\medbreak
\demo {Proof} For each $i$ there is a  UCT sequence
$$
0   \to \ext {K_*(A_i)}.{K_*(B)}  \to  \KKgraded *.{ A_i }.{ B  }
\to  \hom {K_*(A_i)}.{K_*(B)}  \to 0
$$
and these form an inverse system of short exact sequences.
This yields a six term $\invlim$ - $\invlimone $ sequence. However,
$$
\invlimone {\ext {G_i}.H   }   \, =  \,  0     .
$$
by the results of Roos \cite{8}.
Then one obtains   the algebraic isomorphism
$$
\invlimone\gamma _i : \invlimone \,\KKgraded *.{ A_i }.{ B  }
\longrightarrow
  \invlimone \,\hom {K_*(A_i)}.{K_*(B)}
$$
The isomorphism is continuous, by construction,   so
 it suffices to prove that the map $\invlimone\gamma _i $
is an open map. Each map
$$
\gamma _i : \KKgraded *.{A_i}.B  \to \hom {K_*(A_i)}.{K_*(B)}
$$
is open (since the UCT splits topologically)
 and thus the map $\prod \gamma _i $ is open in the
commutative diagram
$$
\CD
\prod\KKgraded *.{A_i}.B     @>{\prod \gamma _i }>>    \prod\hom
{K_*(A_i)}.{K_*(B)}  \\
@VV\Psi V       @VV\Psi ' V     \\
 \prod\KKgraded *.{A_i}.B     @>{\prod \gamma _i }>>  \prod\hom
{K_*(A_i)}.{K_*(B)}   \\
@VV{\pi }V    @VV{\pi '}V    \\
\invlimone\KKgraded *.{A_i}.B    @>{\invlimone\gamma _i}>>
  \invlimone \hom {K_*(A_i)}.{K_*(B)}    \\
@VVV    @VVV   \\
0   @.    0
\endCD
$$
\medbreak
\flushpar The maps $\pi $ and $\pi '$ are continuous and open
quotient maps. This
shows that $ (\invlimone\gamma _i)\pi  $ is open, which implies that
 $\invlimone\gamma _i $ is open. \qed\enddemo

\medbreak

We note that there is another possible topology on
$\KKgraded *.A.B $, namely that obtained by taking the
sequence of subgroups
$$
Ker {\biggl[}  \KKgraded *.A.B   ] \longrightarrow  \KKgraded *.{A_i}.B
{\biggl ]}
$$
and the induced topology on $\usualext $ obtained
by
the associated sequence of subgroups
$$
Ker  {\biggl [}  \usualext    \longrightarrow   \ext {K_*(A_i)}.{K_*(B)}
{\biggl ]}  .
$$
We denote this topology by using $M$ as a subscript.
\medbreak
The following Proposition will be used in the proof of Proposition 4.4.
Recall that
for any topological group $G$,   $G_o $ denotes  the closure of zero,
and $\overline G = G/G_o $, denotes the maximal Hausdorff quotient
 group of $G$.
\medbreak

\proclaim{Proposition 4.3} Suppose given a commutative diagram
of pseudopolonais groups
$$
\CD
0  @>>>   G_o         @>>>      G       @>{\pi _G}>>  \overline G  @>>>  0 \\
@.  @VV{\lambda _o}V    @VV{\lambda}V    @VV{\bar\lambda}V  \\
0  @>>>   H_o         @>>>      H       @>{\pi _H}>>  \overline H  @>>>  0 \\
\endCD
$$
with $\bar\lambda $ an isomorphism of topological groups and
$\lambda _o $ and $\lambda $ algebraic isomorphisms. Then
$\lambda $ and $\lambda _o$ are  isomorphisms of topological
groups. Consequently, any algebraic splitting $\overline G \to G$
is continuous.
\endproclaim

\medbreak
\demo{Proof} The map $\lambda _o$ is an algebraic isomorphism
of topological groups with the indiscrete topology and is
automatically a homeomorphism.  So we
concentrate upon $\lambda $.

Let $U$ be an open neighborhood of the origin in $H$. Then
 $U = U + H_o = \pi _H ^{-1}\pi _H \,U$ is
saturated and hence
$\pi _H (U)$ is open in $\overline H$.
 Then
$$
\lambda ^{-1}U
= \pi _G^{-1}\bar\lambda ^{-1} \pi _H U
= (\bar\lambda \pi _G)^{-1}\pi _H U.
$$
The map $ \bar\lambda \pi _G$ is continuous, $\pi _H C $
is open, and hence $\lambda ^{-1}C$ is open in $G$. Thus
$\lambda $ is continuous.

Reversing the roles of $G$ and
$H$ by symmetry, $\lambda ^{-1} $ is also continuous, and hence
$\lambda $ is an isomorphism of topological groups.
\qed\enddemo

\medbreak

\proclaim{Proposition 4.4} There are  natural isomorphisms
of topological groups
$$
\align
\usualext _M  &\cong  \usualext _J   \\
              &\cong \usualext _{rel}  \\
                &\cong \usualext _I
\endalign
$$
\endproclaim

\medbreak
\demo{Proof} The
first isomorphism
is a consequence of the technique of $KK$-filtration,
since any increasing sequence of subgroups of $K_*(A)$ may
be realized as a sequence of the form $\{\, K_*(A_i)\,\}$.
For the second we argue as follows. The natural map
$$
\iota : \usualext _J  \to \usualext _{rel}
$$
is continuous and induces a commuting diagram
$$
\minCDarrowwidth{.155in}
\CD
0 @>>> \zzero J  @>>> \usualext _J
@>>> \invlim\ext{K_*(A_i)}.{K_*(B)}
@>>> 0   \\
@.  @VV{\iota '}V   @VV{\iota}V   @VV{\iota ''}V   \\
0 @>>> \zzero {rel}  @>>> \usualext _{rel}
 @>>> \invlim\ext{K_*(A_i)}.{K_*(B)}
@>>> 0.
\endCD
$$
We deliberately do not put subscripts on the term
$\invlim\ext {K_*(A_i)}.{K_*(B) }  $
since we know that the topologies on this term are homeomorphic.
Now $\iota '$ is an algebraic and hence a topological isomorphism,
and $\iota ''$ is an isomorphism of topological groups as well.
So we apply Proposition 4.3 to conclude that $\iota $ is
also an isomorphism of topological groups.

A similar argument shows that the natural map
$$
\iota : \usualext _I  \to \usualext _{rel}
$$
is an isomorphism of topological groups. The key point is
that $\iota $ is continuous, and this follows as in (3.4)
\qed\enddemo
\medbreak

 We assume as usual that $A \in \Cal N$, so
that the UCT holds. This takes the form of a natural short exact
sequence  \cite{7}
$$
0 \to
\ext{K_*(A)}.{K_*(B)}  \longrightarrow \KKgraded *.A.B
\to
\hom{K_*(A)}.{K_*(B)}
\to 0
$$
and  $\KKgraded *.A.B $ has a natural structure as a
pseudo\-polonais topological group,
\medbreak

\proclaim {Theorem 4.5}
 Suppose that $A \in \Cal N$. Then each of the
  splittings   of the Universal Coefficient Theorem  constructed in
\cite {9}
  is continuous, and
for each splitting
the resulting (unnatural) algebraic isomorphism
$$
\KKgraded *.A.B  \,\cong\, \hom {K_*(A)}.{K_*(B)}    \oplus
\ext  {K_*(A)}.{K_*(B)} _{rel}
$$
is an isomorphism of pseudo\-polonais  groups.
If $K_*(A) $ is finitely generated then the group  $\KKgraded *.A.B $
is polonais.

\endproclaim
\medbreak
\demo {Proof} The splittings are constructed as follows. Choose
$KK$-equivalences
$$
A \,\approx\, A_0 \,\oplus \, A_1  \qquad\qquad
B \,\approx\, B_0 \,\oplus \, B_1      .
$$
with
 $$
K_i(A_j) \,=\,  K_i(B_j) \,=\, 0  \qquad i \neq j
$$
Then the UCT breaks down into the direct sum of four sequences and
since a $KK$-equivalence is a homeomorphism
it is enough to consider two of these cases, say $A_0 ,\,\,B_0 $
  and    $A_0 ,\,\,B_1 $, and furthermore we may assume that
  each $A_i$ is $KK$-filtered.

First assume that $A = A_0$ and $B = B_0$.
Then  there is only one UCT map which is nontrivial, namely
$$
\gamma : \KKgraded 0.{A }.{B }
\longrightarrow
\hom {K_*(A )}.{K_*(B)} _0\,\cong\, \hom {K_0(A )}.{K_0(B )}      .
 $$
This map is an isomorphism, and it is
an isomorphism of topological groups  by \cite {13, 4.6}.
 So $\Gamma = \gamma ^{-1} $
is a continuous splitting.

In the case $A = A_0$ and $B = B_1$ the map
$$
\delta :  \ext {K_*(A)}.{K_*(B)} _{rel}
  \longrightarrow \KKgraded *.{A}.{B}
$$
is an isomorphism, by the UCT, and it is
continuous. To show that it is an isomorphism we consider the
commuting diagram with exact columns:
\medbreak
$$
\CD
0   @.           0       \\
@VVV        @VVV        \\
\invlimone\hom {K_*(A_i)}.{K_*(B)} @>{\delta '}>>
\invlimone\hom {K_*(A_i)}.{K_*(B)}  \\
@VVV         @VVV       \\
\ext {K_*(A)}.{K_*(B)}    @>{\delta }>>
\ext {K_*(A)}.{K_*(B)}   \\
@VVV            @VVV          \\
\invlim\ext {K_*(A_i)}.{K_*(B)}     @>{\delta ''}>>
\invlim\KKgraded *.{A_i}.{B}  \\
@VVV            @VVV           \\
0        @.          0
\endCD
$$
\medbreak\medbreak\flushpar
The horizontal maps are algebraic isomorphisms. The map $\delta '$
is the identity map of topological groups. The map $\delta ''$
is an isomorphism by Proposition 3.2.
Thus we may apply \cite {13, Th. 6.5} to conclude that
the map $\delta $ is an isomorphism of topological groups.

Finally, if $K_*(A)$ is finitely generated then $Ext $ is countable,
by 1.5, and the result follows.
 \qed\enddemo
\medbreak

{\bf {Remark 4.6}.} We may realize the Jensen example (2.4) as follows.
Choose $C^*$-algebras $A, B \in \Cal N$ with
$$
K_0(A) = \oplus _1^\infty \Bbb Z /2  \qquad   K_0(B) = \Bbb Z/2
\qquad K_1(A) = K_1(B) = 0.
$$
(Such $C^*$-algebras exist and are unique up to $KK$-equivalence
by the proof of the UCT.) Then
$$
\KKgraded 1.A.B \cong \ext{K_0(A)}.{K_0(B)} \cong
\ext{\oplus _1^\infty \Bbb Z /2 }.{\Bbb Z/2}  \cong
\prod _1^\infty \Bbb Z/2
$$
and
applying the analysis of 2.4 we see that
 $\usualext    _{\Bbb Z} $
is not homeomorphic to $\usualext _J$.

\newpage
\beginsection{5. Splitting Obstructions}

In this section we introduce splitting obstructions
$m(A,B)$ and $j(A,B)$ associated to the Milnor and Jensen
sequences respectively and show how they are related.
We demonstrate that if either $K_*(A)$ or $K_*(B)$ is torsionfree
then both sequences split. We also give an example due to
Christensen-Strickland in which the obstructions do not vanish.

We suppose as usual that $A \in \Cal N$ with associated
KK-filtration diagram (1.4).
 Define
$$
m(A,B) \in \ext {  \, \invlim \KKgraded *.{A_i}.B \,    }.{\,\invlimone
\KKgraded *.{A_i}.B \, }
$$
to be the class of the Milnor sequence
$$ 0 \to\invlimone \KKgraded *.{A_i}.B
\to \KKgraded *.A.B  \overset\rho\to\rightarrow
 \invlim \KKgraded *.{A_i}.B \to 0
$$
and let
 $$
j(A,B) \in
\ext{    \invlim\ext{K_*(A_i)}.{K_*(B)}   }.{   \pext {K_*(A)}.{K_*(B)}   }  .
$$
be the class of the Jensen sequence
$$
\minCDarrowwidth {.155in}
\CD
0 @>>> \usualpext  @>>>  \usualext  @>\varphi >> \invlim\ext{K_*(A_i)}.{K_*(B)}
@>>>    0
\endCD
$$
Thus the Milnor sequence splits  iff $m(A,B) = 0$
and similarly the Jensen sequence splits  iff
$j(A,B)= 0$.\footnote{ The maps $\rho $ and $\varphi $ are
continuous and open. For each map an algebraic splitting
is automatically continuous, by Proposition 4.5.
No claim is made for the naturality
of such splittings.}

Note that we have shown that
$$
\invlimone \KKgraded *.{A_i}.B  \,\cong\, \usualpext
$$
and so the subgroups in the two short exact sequences are
canonically isomorphic. We make this identification in
 the following proposition.

 Here is how the two splitting obstructions  are related.

\medbreak

\proclaim{Proposition 5.1} With the notation above, and up to
the Jensen
isomorphism,
$$
j(A,B) = ( \invlim\delta _i)^* m(A,B)   .
$$
\endproclaim
\medbreak
\demo{Proof} In light of the Jensen isomorphism it suffices to
demonstrate that the right square of the commutative diagram
$$
\minCDarrowwidth {.155in}
\CD
0 @>>> \usualpext  @>\psi >>  \usualext  @>\varphi >> \invlim\ext{K_*(A_i)}.{K_*(B)}
@>>>    0  \\
@.   @VV{\cong }V    @VV\delta V   @VV{\invlim \delta _i }V   \\
 0 @>>>   \invlimone \KKgraded *.{A_i}.B
@>\sigma >>   \KKgraded *.A.B  @>\rho >>
\invlim  \KKgraded *.{A_i}.B @>>>     0
 \endCD
$$
is a pullback diagram.

So suppose that $x \in \KKgraded *.A.B $,
$y \in  \invlim\ext{K_*(A_i)}.{K_*(B)} $,
and that
$$
\rho (x) = (\invlim \delta _i )(y).
$$
 Then
$$
\align
\gamma (x)  &= ( \invlim \gamma _i)\rho (x)  \\
&= (\invlim \gamma _i )( \invlim \delta _i )(y)   \\
&= \invlim (\gamma _i\delta _i)(y)  \\
&= 0
\endalign
$$
since $\gamma _i \delta _i  = 0$ for each $i$, and hence
$$
x \in Ker( \gamma ) = Im (\delta ).
$$
 As $\delta $ is mono,
there is a unique choice $z \in \usualext $ with $\delta (z) = x$.
Then
$$
\align
(\invlim \delta _i )( y - \varphi (z) ) &= \rho (x) -
(\invlim \delta _i)\varphi (z) \\
&= \rho (x) - \rho \delta (z)  \\
&= \rho (x) - \rho (x) = 0
\endalign
$$
and since $ \invlim \delta _i $ is mono, this implies that $y = \varphi (z)$.
Thus the right square is a pullback and the proof is complete.
 \qed\enddemo
\medbreak
 Here is another result
along similar lines.

\medbreak

\proclaim{Theorem 5.2}
\roster
\item
If $K_*(A)$ is a direct sum of cyclic groups or if $K_*(B)$
is algebraically compact then $\usualpext = 0$  and
  both obstructions vanish.
\medbreak\item
If $K_*(A)$ or $K_*(B)$ is  torsionfree
then
 $\usualpext $ is divisible and both obstructions vanish.
\endroster
\endproclaim

\medbreak
\demo{Proof}
These are both purely algebraic results. Let $G = K_*(A)$ and
$H = K_*(B)$.
If $G$ is a direct sum of cyclic groups or if $H$
is algebraically compact  then the group
$\pext G.H = 0$ by classical results
(cf. \cite{16}, 5.4).

If $G$ is torsionfree then for any abelian group $H$ the
group $\ext G.H $ is divisible. Then
$$
\pext G.H = \cap _n n\ext G.H     =\ext G.H
$$
and so $\pext G.H $ is divisible.
If $H$ is torsionfree then for any $G$ the group $\pext G.H $
is the maximal divisible subgroup of $\ext G.H $  by
\cite{16, 8.5}.  This completes the proof.
\qed\enddemo

\medbreak

To look for cases where the sequences do not split, then,
one must have torsion (for the same prime) in both
$K_*(A)$ and $K_*(B)$.

\medbreak

{\bf Example 5.3.} (Christensen-Strickland \cite{3})
\medbreak

Here is an example where both invariants are non-vanishing.
 Fix a prime $p$.
Now take $A,\, B \in \Cal N$ with
$$
K_0(A) = \Bbb Z(p^{\infty})  \qquad K_1(A) = 0
$$
$$
K_0(B) = 0  \qquad   K_1(B) =  \underset{n}\to\oplus\,  \, \Bbb Z/p^n
$$
Then $\hom {K_*(A)}.{K_*(B)} = 0 $ and
the Milnor and Jensen sequences both reduce to the short exact
sequence
$$
0 \to \pext {K_0(A)}.{K_1(B)}
\to \KKgraded 0.A.B
\to  \invlim   K_1(B)/p^n
\to 0.
$$
Making algebraic identifications, this sequence is isomorphic
to the sequence
$$
0 \to  \pext {\Bbb Z(p^{\infty}) }.{ \underset{n}\to\oplus\,  \Bbb Z/p^n}
\to \ext {\Bbb Z(p^{\infty}) }.{ \underset{n}\to\oplus\, \Bbb Z/p^n}
\to  \widehat {\underset{n}\to\oplus\,  \, \Bbb Z/p^n}
\to 0
$$
where
$\widehat {\underset{n}\to\oplus\,  \, \Bbb Z/p^n}$
denotes the $p$-adic completion of
${\underset{n}\to\oplus\,  \, \Bbb Z/p^n}$.
Christensen and Strickland \cite{3,  6.6, 6.7} demonstrate that
this sequence does not split. Thus $j(A,B) \neq 0$\footnote{In fact
they show that the element $j(A,B)$ has infinite order.} and
by Prop. 5.1, $m(A,B) \neq 0$ as well.

 \ref\key   {\bf  }
\by
\paper
\jour
 \vol
\yr
\pages
\endref
\medbreak

\end
\newpage
\Refs
\widestnumber\key{XXXXXXX}






 \ref\key   {\bf 1}
\by   B. Blackadar
\book     K-Theory for Operator Algebras
\publ  Math. Sci. Res. Inst. No. 5, 2nd. Ed., Cambridge U. Press
\publaddr New York
\yr     1998
\endref\medbreak




 \ref\key   {\bf 2}
\by   H. Cartan and S. Eilenberg
\book     Homological Algebra
\publ Princeton University Press
\publaddr Princeton
\yr     1956
\endref\medbreak



 \ref\key   {\bf 3}
\by       J. D. Christensen and N. P. Strickland
\paper   Phantom maps and homology theories
 \jour Topology
\vol 37
\yr  1988
\pages  339-364
\endref
\medbreak


\ref\key   {\bf 4}
\by       M. Dadarlat
\paper   Approximate unitary equivalence and the topology
of $Ext (A,B) $
 \jour
\endref
\medbreak


 \ref\key   {\bf 5}
\by    L\'aszl\'o Fuchs
\book      Infinite  Abelian  Groups,
       {\rm Pure and Applied Mathematics No. 36}
\publ Academic Press
\publaddr New York
\vol 1
\yr     1970
\pages  290
\endref\medbreak



  \ref\key   {\bf 6}
\by   C. U. Jensen
\book     Les Foncteurs D\'eriv\'es de $\invlim $ et leur
Applications en Th\'eorie des Modules,
       {\rm Lecture Notes in Mathematics}
\vol 254
\publ Springer Verlag
\publaddr New York
\yr     1972
\endref\medbreak

 \ref\key   {\bf 7}
\by    S. MacLane
\book      Homology
\publ Springer-Verlag
\publaddr New York
\yr   3rd ed.   1975
\endref\medbreak




  \ref\key   {\bf   8}
\by   J. E. Roos
\paper   Sur les d\'eriv\'es de $\invlim $
\jour   C. R. Acad. Sci. Paris
\vol   252
\yr      1961
\pages      3702-3704
\endref
\medbreak




 \ref\key   {\bf 9}
\by     J. Rosenberg and C. Schochet
\paper    The K\"unneth theorem and the universal coefficient theorem for
             Kasparov's generalized K-functor
\jour     Duke Math. J.
 \vol     55
\yr     1987
\pages     431-474
\endref
\medbreak



  \ref\key   {\bf   10}
\by   N. Salinas
\paper   Relative quasidiagonality and $KK$-theory
\jour   Houston J. Math.
\vol   18
\yr      1992
\pages      97-116
\endref
\medbreak






  \ref\key   {\bf  11}
\by       C. Schochet
\paper  The UCT, the Milnor sequence, and a canonical decomposition
of the Kasparov groups
\jour K-Theory
 \vol   10
 \yr   1996
 \pages   49-72
\endref
\medbreak

\ref\key   {\bf 12}
\by       C. Schochet
\paper  Correction to: The UCT, the Milnor sequence, and a canonical
decomposition
of the Kasparov groups
\jour K-Theory
 \vol   14
 \yr   1998
 \pages   197-199
\endref
\medbreak



  \ref\key   {\bf 13}
\by       C. Schochet
\paper  The fine structure
of the Kasparov groups I: continuity of the $KK$-pairing
 \jour  J. Functional Analysis, to appear
\endref
\medbreak


  \ref \key  {\bf 14}
\by       C. Schochet
\paper  The fine structure
of the Kasparov groups III: relative quasidiagonality
 \jour  submitted
\endref
\medbreak


  \ref\key   {\bf 15}
\by       C. Schochet
\paper  The topological Snake Lemma and corona algebras
 \jour  New York J. Math.
\vol  5
\yr  1999
\pages  131-137
\endref
\medbreak

 \ref\key   {\bf 16}
\by       C. Schochet
\paper   A Pext Primer: Pure extensions and $\invlimone $
for infinite abelian groups
 \jour  submitted
\endref
\medbreak







\endRefs
